\documentclass[a4paper]{mathscan}
 
       \newtheorem{thm}{Theorem}[section] 
       \newtheorem{pro}[thm]{Proposition} 
           
       \newtheorem{lem}[thm]{Lemma}        

       \theoremstyle{definition} 
         
       \newtheorem{rem}[thm]{Remark}   
       \newtheorem{defn}[thm]{Definition}  
       \newtheorem{exam}[thm]{Example}
       \newtheorem{hyp}[thm]{Hypothesis}

\newcommand{\supp}{\operatorname{supp}}

\newcommand{\rn}{{\mathbb R}^n}
\newcommand{\rnp}{{\mathbb R}^n_+}

\newcommand{\rnpm}{\mathbb R^n_\pm}

\newcommand{\comega}{\overline\Omega }
\newcommand{\ang}[1]{\langle {#1} \rangle}
\newcommand{\Op}{\operatorname{Op}}
\newcommand{\simto}{\overset\sim\rightarrow}
\newcommand{\ol}{\overline}
\newcommand{\N}{\mathbb N}
\newcommand{\R}{\mathbb R}
\newcommand{\C}{\mathbb C}

\newcommand{\F}{\mathcal F}

\begin{document}

\title[Weyl asymptotics for fractional-order operators]
{Weyl asymptotics for fractional-order  Dirichlet realizations in
nonsmooth cases}

\author {Gerd Grubb}

\address
{Dept. of Mathematical Sciences,\\ Copenhagen University,\\
Universitetsparken 5,\\ DK-2100 Copenhagen, Denmark.\\
E-mail {\tt grubb\@math.ku.dk}}
\begin{abstract}
Let $P$ be a symmetric $2a$-order classical strongly elliptic pseudodifferential
operator with \emph{even} symbol $p(x,\xi )$ on $\rn$ ($0<a<1$), for example a
perturbation of $(-\Delta )^a$. Let $\Omega \subset\rn$ be  bounded, and let $P_D$ be the Dirichlet realization in
$L_2(\Omega )$ defined under the exterior condition $u=0$ in
$\rn\setminus\Omega $. When $p(x,\xi )$ and $\Omega $ are $C^\infty $, it is known
that the eigenvalues $\lambda _j$ (ordered in a nondecreasing sequence
for $j\to\infty $) satisfy a Weyl asymptotic formula
\begin{equation*}
\lambda 
_j(P_{D})=C(P,\Omega )j^{2a/n}+o(j^{2a/n})\text{
for }j\to \infty,
\end{equation*}
with $C(P,\Omega )$ determined from the principal symbol of $P$. We
now show that this result is valid for more general operators with a
possibly nonsmooth $x$-dependence, over Lipschitz domains, and that it extends to $\tilde
P=P+P'+P''$, where $P'$ is an operator of order $<\min\{2a,
a+\frac12\}$ with certain mapping properties, and $P''$ is bounded in $L_2(\Omega )$ (e.g.\
 $P''=V(x)\in L_\infty (\Omega )$). Also the regularity of
eigenfunctions of $P_D$ is discussed.
\end{abstract}

\subjclass{35S15; 47G30; 35J25; 60G51} 

\dedicatory{Dedicated to Vsevolod A.\ Solonnikov, on the occasion of his
  ninetieth birthday}

\maketitle

\section { Introduction}\label{sec1}

Consider a $2a$-order classical strongly elliptic  pseudodifferential operator $P$ on
$\rn$ ($0<a<1$) with a  symbol $p(x,\xi )$ with \emph{even} parity in
$\xi $ (cf.\ \eqref{eq:2.5}), having $C^\tau
$-dependence on $x$ for some $\tau >2a$. A standard example is the fractional
Laplacian $(-\Delta )^a$, which is $x$-independent, but our methods
also allow $x$-dependent symbols. Let $\Omega \subset\rn$ be bounded
open with $C^{1+\tau }$-boundary. Denote by $P_D$ the $L_2$-Dirichlet
realization of $P$ with domain $D(P_{D})= \{u\in \dot H^a(\comega)\mid (Pu)|_\Omega \in
L_2(\Omega )\}$. Its spectrum $\Sigma $ is discrete and contained in a convex
sectorial region in $\C$ opening to the right.

Section 2 contains preliminary material.

In Section 3, we recall the results from \cite{G23} on the regularity of the eigenfunctions $u_\lambda $:
\begin{equation}\label{eq:1.1}
u_\lambda \in   d^aC^{t }(\comega)\text{ for
}t=\min\{2a,\tau -a\}-\varepsilon ;\text{ small }\varepsilon >0,
\end{equation}
where $d(x)=\operatorname{dist}(x,\partial\Omega )$. Also the
$L_q$ Dirichlet realizations are considered for $1<q<\infty $; they have the
same spectrum and eigenfunctions as for $q=2$.
Slightly sharper results than  \eqref{eq:1.1} hold in $C^\infty
$-cases \cite{G15b}.

In Sections 4 and 5, we study the asymptotic behavior of the eigenvalues, restricting
the attention to 
selfadjoint operators. The main purpose is to
show that the eigenvalues $\lambda _j(P_D)$, ordered nondecreasingly and repeated
according to multiplicity, satisfy
\begin{align}\label{eq:1.2}
\lambda 
_j(P_{D})&=C(P,\Omega )j^{2a/n}+o(j^{2a/n})\text{
for }j\to \infty ,\text{ where}\\ \nonumber
C(P,\Omega )&=\Bigl(\tfrac 1{n(2\pi )^n}\int_{\Omega }\int_{|\xi
              |=1}|p_{0}(x,\xi )|^{-n/2a}\,d\omega (\xi
              )dx\Bigr)^{-2a/n};
\end{align}
 with $p_0$ denoting the principal symbol of $P$. 

This was shown for $(-\Delta )^a$ already by Blumenthal and Getoor in
\cite{BG59} when $\Omega $ has bounded volume. See also e.g.\ Chen and Song
\cite{CS05}, Geisinger \cite{Ge14}, Dyda, Kuznetsov and Kwasnicki \cite{DKK17}. Also Frank and others contributed to
the question, see the survey Frank \cite{F17} for an account of
related eigenvalue estimates (such as the behavior of sums of
eigenvalues), including fractional Schr\"odinger operators $(-\Delta
)^a+V(x)$. Two-term asymptotics were put forward by Ivrii \cite{I16}.
In
\cite{G15b}, we established \eqref{eq:1.2} for $x$-dependent operators $P$ with smooth symbols over smooth
regions $\Omega $, also giving estimates of $s$-numbers (singular
values) in nonsymmetric cases. 

Our present aim is to allow nonsmooth regions $\Omega $ as well as
more general operators $P$
with possibly nonsmooth $x$-dependence. For the main results we focus on selfadjoint cases,
where there are convenient  perturbation tools. Nonsmooth symbols
$p(x,\xi )$ are treated in Section 4, and estimates when $\Omega $ is
Lipschitz are obtained in Section 5.

Besides $P$, we consider suitable lower-order perturbations
\begin{equation}\label{eq:1.3}
\tilde P=P+P'+P'',
\end{equation}
where $P'$ is of
order $<\min\{2a,a+\frac12\}$ with certain mapping properties in
Sobolev spaces, and $P''$ is bounded in $L_2(\Omega )$
(e.g.\ the multiplication by
a  potential $V(x)\in L_\infty (\Omega )$). Then 
the Dirichlet realization $\tilde P_D$  is  shown to have the same
domain $H^{a(2a)}(\comega)$  as $P_D$ 
(Theorem \ref{Theorem 4.6} $1^\circ$), and  
\eqref{eq:1.2} is shown to hold for $\tilde P_D$ 
 with the same constant $C(P,\Omega )$ when $\Omega $ is Lipschitz
 (Theorem \ref{Theorem 5.3}).

\section {Preliminaries}\label{sec2} 


The setting is the same as in \cite{G23}, so we shall just rapidly
repeat the part of the notation that will be needed here.

The
space $C^k(\rn)\equiv C^k_b(\rn)$ consists  of $k$-times differentiable
functions with bounded norms $\|u\|_{C^k}=\sup_{|\alpha |\le k,x\in\rn}|D^\alpha
u(x)|$ ($k\in{\mathbb N}_0=\{0,1,2,\dots\}$). The spaces
$C^{k,\sigma } (\rn)$, $k\in \N_0$ and $0<\sigma \le 1$,
 consist of
function  with bounded norms
$\|u\|_{C^{k,\sigma } }=\|u\|_{C^k}+\sup_{|\alpha |= k,x\ne y}|D^\alpha
u(x)-D^\alpha u(y)|/|x-y|^\sigma $. For $\sigma <1$, they are the H\"older
spaces, also denoted $C^{k+\sigma } (\rn)$.
For $\sigma =1$, they are Lipschitz spaces. Occasionally, we refer to the
scale of H\"older-Zygmund spaces $C_*^s(\rn)$, $s\in\R$,
which coincide with the H\"older spaces for $s>0, s\notin \N$, and is
preserved under interpolation.
There are similar spaces over subsets
of $\rn$.

The halfspaces $\rnpm$ are defined by
 $\rnpm=\{x\in
{\mathbb R}^n\mid x_n\gtrless 0\}$, with points denoted  $x=(x',x_n)$,
$x'=(x_1,\dots, x_{n-1})$. For a given real function $\zeta \in C^{t}
(\R^{n-1})$ (some $t \ge 0$), we define the curved halfspace
$\rn_\zeta  $ by $
\rn_\zeta = \{x\in\R^n\mid x_n>\zeta (x')\}$;
it is
 a $C^{t }$-domain. 
By a bounded $C^{t }$-domain $\Omega $  we mean the following:
 $\Omega \subset\rn$ is open and bounded, and every boundary point
$x_0$ has an open neighborhood $U$ such that, after a translation of
$x_0$ to $0$ and a suitable rotation, $U\cap \Omega $ equals $U\cap \rn_\zeta 
$ for a function $\zeta  \in C^{t }(\R^{n-1})$ with $\zeta
(0)=0$. There are similar definitions for $C^{k,1}$-spaces.

Restriction from $\R^n$ to $\rnpm$ (or from
${\mathbb R}^n$ to $\Omega $ resp.\ $\complement\comega= \rn \setminus \comega$) is denoted $r^\pm$,
 extension by zero from $\rnpm$ to $\R^n$ (or from $\Omega $ resp.\
 $\complement\comega$ to ${\mathbb R}^n$) is denoted $e^\pm$.

When $\Omega $ is a $C^{1+\tau }$-domain with $\tau >0$, we denote by $d(x)$ 
a function  that is $C^{1+\tau } $ on $\comega$,
positive on $\Omega $ and vanishes only to the first order on
$\partial\Omega $ (i.e., $d(x)=0$ and $\nabla d(x)\neq 0$ for $x\in
\partial\Omega$). It is equivalent to the distance
 $d_0(x)=\operatorname{dist}(x,\partial\Omega )$ near $\partial\Omega
 $. (More details in \cite{AG23,G23}.)

The Bessel-potential spaces
$H^s_q({\mathbb R}^n)$ are defined
for $s\in{\mathbb R}$, $1<q<\infty $, by 
\begin{equation}\label{eq:2.1}
H_q^s(\R^n)=\{u\in \mathcal S'({\mathbb R}^n)\mid \F^{-1}(\ang{\xi }^s\hat u)\in
L_q(\R^n)\},
\end{equation}
where $\ang\xi=(|\xi |^2+1)^{\frac12}$, and $\mathcal F$ is the Fourier transform  $\hat
u(\xi )=(\mathcal F
u)(\xi )= \int_{{\mathbb R}^n}e^{-ix\cdot \xi }u(x)\, dx$, with inverse $(\mathcal F^{-1}
v)(x)= (2\pi )^{-1}\int_{{\mathbb R}^n}e^{+ix\cdot \xi }v(\xi )\, d\xi
$.  For $s\in \N_0$, the spaces are also denoted
$W^{s,q}(\rn)$. When  $q=2$, they are the standard $L_2$ Sobolev spaces,
where the index $2$ is usually omitted.

Along with the spaces $H^s_q({\mathbb R}^n)$ defined in \eqref{eq:2.1}, there are
the two scales of spaces associated with $\Omega $ for $s\in{\mathbb R}$:
\begin{align} \nonumber
\ol H_q^{s}(\Omega)=&\{u\in \mathcal D'(\Omega )\mid u=r^+U \text{ for some }U\in
                     H_q^{s}(\R^n)\},\\ \label{eq:2.2}
  &\text{ the \emph{restricted} space},\\ \nonumber
\dot H_q^{s}(\comega)=&\{u\in H_q^{s}({\mathbb R}^n)\mid \supp u\subset
\comega \},\text{ the \emph{supported} space;}
  \end{align}
here $\supp u$ denotes the support of $u$ (the complement
 of the largest open set where $u=0$). 

A \emph{pseudodifferential operator} ($\psi $do) $P$ on ${\mathbb R}^n$ is
defined from a function $p(x,\xi )$  on ${\mathbb
R}^n\times{\mathbb R}^n$, called the \emph{symbol},  by 
\begin{equation}\label{eq:2.3}
Pu=\Op (p(x,\xi ))u 
=(2\pi )^{-n}\int_{\rn} e^{ix\cdot\xi
}p(x,\xi )\hat u(\xi)\, d\xi =\mathcal F^{-1}_{\xi \to x}(p(x,\xi )\F u(\xi
)),
\end{equation}
using the Fourier transform $\F$. An introduction to $\psi $do's is
given e.g.\ in \cite{G09}, Ch.\ 7--8, and a description with further references
and an inclusion of results for operators with nonsmooth symbols can
be found in
\cite{AG23}. 

The space $S^m_{1,0}({\mathbb R}^n\times{\mathbb R}^n)$ of symbols
$p$ of order $m\in{\mathbb R}$ consists of the complex
$C^\infty $-functions $p(x,\xi )$
such that $\partial_x^\beta \partial_\xi ^\alpha p(x,\xi
)$ is $O(\ang\xi ^{m-|\alpha |})$ for all $\alpha ,\beta $, for some
$m\in{\mathbb R}$, with global estimates in $x\in{\mathbb R}^n$.
$P$ is then of order $m$; it  maps $H^s_q({\mathbb R}^n)$ continuously into
$H^{s-m}_q ({\mathbb R}^n)$ for all $s\in{\mathbb R}$.

$P $ with symbol $p\in  S^m_{1,0}({\mathbb R}^n\times{\mathbb R}^n)$ is said to be \emph{classical} when 
$p$  
has an asymptotic expansion $p(x,\xi )\sim \sum_{j\in{\mathbb
N}_0}p_j(x,\xi )$ with $p_j$ \emph{homogeneous in} $\xi $ of degree $m-j$ for
all $|\xi |\ge 1$ and $j\in\N_0$, such that
\begin{equation}\label{eq:2.4}
\partial_x^\beta \partial_\xi ^\alpha \bigl(p(x,\xi )-
{\sum}_{j<J}p_j(x,\xi )\bigr) \text{ is }O(\ang\xi ^{m-|\alpha |-J})\text{ for
all }\alpha ,\beta \in{\mathbb N}_0^n, J\in{\mathbb N}_0.  
\end{equation}
The subspace of classical symbols is denoted $S^m({\mathbb R}^n\times{\mathbb R}^n)$.

A classical symbol $p(x,\xi )$ of order $m\in\R$ (and the associated operator $P$) is said to be   \emph{strongly elliptic}
when $\operatorname{Re}p_0(x,\xi )\ge c|\xi |^m $ for $|\xi |\ge 1$,
with $c>0$.
Moreover, a classical $\psi $do $P=\Op(p(x,\xi ))$ of
order $m$ is said
to be 
\emph{even}, when the terms in the symbol expansion $p\sim\sum_{j\in\N_0}p_j$ satisfy
\begin{equation}\label{eq:2.5}
p_j(x,-\xi )=(-1)^jp_j(x,\xi )\quad \text{ for all }x\in \rn,|\xi |\ge
1, j\in \mathbb N_0.
\end{equation}

There are also defined spaces of symbols with finite smoothness in
$x$. In the case of  $C^\tau $-smoothness ($\tau >0$), the symbol
spaces are denoted $C^\tau S^m_{1,0}(\R^{n}\times \rn)$, $C^\tau
S^m(\R^{n}\times \rn)$, see details  e.g.\ in
\cite{AG23,G23}. Also here, strongly elliptic symbols and even
symbols are well-defined subclasses.
When $p\in C^\tau
S^m_{1,0}(\R^{n}\times \rn)$, it defines a continuous operator 
\begin{equation}\label{eq:2.6}
\Op(p)\colon H^{s+m}_q (\R^n)\to H^s_q(\R^n)\text{  for }|s|<\tau .
\end{equation}

There  are two families of pseudodifferential operators that play a special
role; the ``order-reducing operators'' $\Lambda _\pm^{(t)}$ defined for bounded $C^\infty
$-domains $\Omega $ in \cite{G15a}. They are of order $t$, and have the
homeomorphism properties:
\begin{equation}\label{eq:2.7}
\Lambda  ^{(t) }_+\colon \dot H_q^s(\comega )\simto
\dot H_q^{s- t }(\comega),\quad
r^+\Lambda  ^{(t) }_{-}e^+\colon \ol H_q^s(\Omega  )\simto
\ol H_q^{s- t } (\Omega  ), \text{ all }s\in\R.
\end{equation}
Here the operator $r^+\Lambda  ^{(t) }_{-}e^+$, also denoted $\Lambda
^{(t) }_{-,+}$, is the adjoint of $\Lambda _+^{(t)}$.
The operators are modeled after the operators  $\Xi
^t_\pm=\operatorname{Op}\bigl((\ang{\xi '}\pm i\xi _n)^t\bigr)$ that act in this way relative to
$\rnp$. (For
$C^{1+\tau }$-domains we have not made the effort to introduce global versions of $\Lambda
_\pm^{(t)}$, but obtain the relevant results by using versions of $\Xi _\pm^t$ in
local coordinates at the boundary.) They play a role in the study of
our special $2a$-order operators $P$, because they  enter as factors
in the analysis of the  homogeneous Dirichlet problem \eqref{eq:3.2}
for $P$. In fact, the solution space for the problem with right-hand
side $f$ in $\ol
H_q^{s}(\Omega )$
is the so-called  \emph{$a$-transmission space}:  
\begin{align}\label{eq:2.8}
H_q^{{a} (s+2a)}(\comega)&=\Lambda  _+^{-(a) }e^+\ol H_q^{s+a}(\Omega
                           )\text{ for }s>-a-1/q',\\ \nonumber
  H_q^{a(s+2a)}(\comega)&=\dot H_q^{s+2a}(\comega)\text{ if }s<-a+1/q,
\end{align}
defined via local coordinates when $\Omega $ is nonsmooth ($\frac1q + \frac1{q'}=1$). This regularity
result is formulated below as Theorem \ref{Theorem 3.3}. A pedestrian introduction
to $a$-transmission spaces can be found in \cite{G22}.

The $a$-transmission spaces are  defined similarly for the scale of
H\"older-Zygmund spaces: $C_*^{{a} (s+2a)}(\comega)=\Lambda  _+^{-(a) }e^+\ol C_*^{s+a}(\Omega )$.

\section {On eigenfunctions and their regularity}\label{sec3}

Our basic hypothesis is:

\begin{hyp}\label{Hypothesis 3.1}
There are given constants $a,\tau ,q$ with $0<a<1$, $\tau >2a$, and 
 $1<q<\infty $.    $P$ is a
 classical, strongly elliptic $\psi $do of order $2a$, with  even symbol in
$C^\tau S^{2a}(\rn\times\rn)$.
\end{hyp}

In much of the paper, we moreover assume that $P$ is symmetric:

\begin{hyp}\label{Hypothesis 3.2} 
The conditions in Hypothesis
\ref{Hypothesis 3.1} hold, and in addition, $P$ equals its formal adjoint $P^*$
(in short: is symmetric).
\end{hyp}

The last hypothesis will be satisfied  e.g.\ if $P$ has real symbol
independent of $x$, if $P$ is a power $L^a$ of a symmetric strongly
elliptic second-order differential operator $L$, or if $P$ is obtained as the ``real part''
$P=\frac12(P_1+P_1^*)$ of a suitable nonsymmetric operator $P_1$.

Recall the regularity theorem for the homogeneous
Dirichlet problem for $P$:

\begin{thm}\label{Theorem 3.3} Assume Hypothesis \ref{Hypothesis 3.1}
and let $\Omega $ be a bounded $C^{1+\tau }$-domain in $\rn$; let
$s$ satisfy  $-a\le s<\tau -2a$.

$1^\circ$ Forward mapping property:
$r^+P$ maps continuously
\begin{equation}\label{eq:3.1}
   r^+ P\colon H_q^{a (s+2a)}(\comega)\to \ol H_q^{s}(\Omega ).
\end{equation}

$2^\circ$ Regularity (backward mapping property):  When
 $u\in
 \dot{H}^a_q(\ol{\Omega})$ solves the homogeneous Dirichlet problem
 \begin{equation}\label{eq:3.2}
   Pu = f \text { in }\Omega ,\quad \supp u\subset \comega,
   \end{equation}
  for some $f\in \ol{H}^s_q(\Omega)$, then $u\in
  H^{a(s+2a)}_q(\ol\Omega)$. Here if $q\ne 2$ and $\tau <\infty $, we assume that $s\ge 0$.
\end{thm}

\begin{proof} The statements were shown in the smooth case ($\tau
=\infty $) for all $1<q<\infty $ in  \cite{G15a} Th.\ 4.2, 4.4 (in fact allowing $s$ to take
values down to 
$>-a-\frac1{q'}$,   and $u\in \dot
H^\sigma _q(\comega)$ with $\sigma >a-\frac1{q'}$ in $2^\circ$). The statements were extended to the nonsmooth cases in
\cite{AG23} Cor.\ 6.2, Th.\ 6.9, with $s$ assumed $\ge 0$ for $2^\circ$. When $q=2$, the interval
$0>s\ge -a$
can be included in $2^\circ$ by use of the variational definition of the Dirichlet
realization, recalled below and studied in more detail in \cite{G23}
Sect.\ 4;
here it is found that $r^+P+b\colon \dot H^a(\comega)\to \ol
H^{-a}(\Omega )$ is a homeomorphism for  suitable
constants $b $. By interpolation between this and the
homeomorphism   $r^+P+b\colon  H^{a(2a)}(\comega)\simto L_2(\Omega )$, we find
$r^+P+b\colon H^{a(2a+s)}(\comega)\simto \ol
H^{s}(\Omega )$ when $0>s>-a$. Then for $-a\le s\le 0$, $u\in
 \dot{H}^a_q(\ol{\Omega})$ and  $r^+Pu\in \ol H^{s}(\comega)$ (hence
 $(r^+P+b)u\in \ol H^{s}(\comega)$) imply $u\in H^{a(2a+s)}(\comega)$.

\end{proof}

Recall from \cite{G23} Sect.\ 4, that under Hypothesis \ref{Hypothesis
  3.1}, $P$
satisfies a G\aa{}rding inequality (with $c>0$, $\beta \in\R$)
\begin{equation}\label{eq:3.3}
\operatorname{Re}(Pu,u)\ge c_0\|u\|^2_{\dot
H^a(\rn)}-\beta \|u\|^2_{L_2(\rn )}, \text{ for }u\in C_0^\infty (\rn),
\end{equation}
and the sesquilinear form $s(u,v)$ obtained by closure in $\dot
H^a(\comega)$ of
\begin{equation}\label{eq:3.4}
s(u,v)=\int_\Omega Pu\,\bar v
\,dx\text{ on }C_0^\infty (\Omega ),
\end{equation}
defines by variational theory (as in e.g.\
\cite{G09} Sect.\ 12.4, after Lions-Magenes \cite{LM68}) the
 $L_2$-Dirichlet 
realization $P_{D,2}$, i.e., the operator in $L_2(\Omega )$
acting like $r^+P$ with domain
\begin{equation}\label{eq:3.5}
D(P_{D,2})= \{u\in \dot H^a(\comega)\mid r^+Pu\in
L_2(\Omega )\}. 
\end{equation}
 Its
spectrum and numerical range is contained in a sectorial region with
$c_1\ge 0$:
\begin{equation}\label{eq:3.6}
M=\{\lambda \in \C\mid \operatorname{Re}\lambda \ge c_0-\beta ,
|\operatorname{Im}\lambda |\le c_1(\operatorname{Re}\lambda +\beta )
\};
\end{equation}
the spectrum is discrete since $\dot H^a(\comega)$ is compactly
injected in $L_2(\Omega )$.
In fact, there holds:

\begin{pro}\label{Proposition 3.4} \cite{G23} Assume Hypothesis \ref{Hypothesis 3.1} and
let $\Omega $ be  bounded and $C^{1+\tau }$. The domain
$D(P_{D,2})$ equals $H^{a(2a)}(\comega)$,  
and the spectrum $\Sigma $ of $P_{D,2}$ is discrete.

For $\lambda \notin \Sigma $, $P_{D,2}-\lambda I $ is a homeomorphism of
$H^{a(2a)}(\comega)$ onto $L_2(\Omega )$.

For $\lambda \in\Sigma $,  $P_{D,2}-\lambda I$ defines a Fredholm
operator with index zero from $H^{a(2a)}(\comega)$ to $L_2(\Omega )$.
The kernel is denoted $N_\lambda $, and there is a cokernel
$N'_{\overline\lambda }$, where $N'_{\overline\lambda }$
denotes the kernel of $P_{D,2}^*-\overline\lambda $. Here 
\begin{equation}\label{eq:3.7}
\operatorname{dim}N_\lambda = \operatorname{dim}N'_{\overline\lambda
}.
\end{equation}
\end{pro}

From \eqref{eq:2.8} with $q=2$, $s=0$, we have that $H^{a(2a)}(\comega)=\Lambda
_+^{-(a)}e^+\ol H^{a}(\Omega )$ (in a local sense if $\Omega $ is
nonsmooth). Here $H^{a(2a)}(\comega)=\dot H^{2a}(\comega )$ if
$a<\frac12$, and $H^{a(2a)}(\comega)\subset\dot
H^{a+\frac12-\varepsilon }(\comega )$ if $a\ge \frac12$ (more
detailed information in \cite{G19,G23}).

We can also define the $L_q$-realization $P_{D,q}$ for general $q$; it acts  as $r^+P$
in $L_q(\Omega )$
with domain
\begin{equation}\label{eq:3.8}
D(P_{D,q})=\{u\in \dot H_q^a(\comega)\mid r^+Pu\in L_q(\Omega
)\}=H_q^{a(2a)}(\comega).
\end{equation}
For this we showed in \cite{G23} Sect.\ 4, that the adjoint of
$P_{D,q}$ is the analogous operator $(P^*)_{D,q'}$ defined from $P^*$
in $L_{q'}(\Omega )$. The spectra and
eigenfunctions satisfy:

\begin{thm}\label{Theorem 3.5} \cite{G23} Assume Hypothesis \ref{Hypothesis 3.1} and let
$\Omega $ be bounded and $C^{1+\tau }$.

$1^\circ$ For any $q\in \,]1,\infty [\,$, the spectrum of $P_{D,q}$
equals $\Sigma $, and when $\lambda \in\Sigma $, the eigenspace equals
$N_\lambda $, and $N'_{\overline \lambda }$ is a cokernel of $P_{D,q}-\lambda $.

$2^\circ$ If $0$ is an eigenvalue, an associated eigenfunction $u_0$ satisfies
\begin{equation}\label{eq:3.9}
u_0\in C_*^{a(\tau -\varepsilon )}(\comega)\subset d^aC^{\tau
-a-\varepsilon }(\comega), 
\end{equation}
for small $\varepsilon >0$ (with  $\tau -a-\varepsilon >0$).

For a nonzero $\lambda \in\Sigma $, the eigenfunctions $u_\lambda $ satisfy
\begin{equation}\label{eq:3.10}
u_\lambda \in  C_*^{a(t+a )}(\comega)\subset d^aC^{t }(\comega)\text{ for
}t=\min\{2a,\tau -a\}-\varepsilon ;\text{ small }\varepsilon >0.
\end{equation}

$3^\circ$ If $\tau =\infty $ or  $P$ is symmetric, the eigenfunctions of $(P^*)_{D,q'}$
have the same behavior as in $2^\circ$. In nonsymmetric cases, when
$\tau $ is finite, they
satisfy at least
\begin{equation}\label{eq:3.11}
u_{\overline\lambda}\ \in  C_*^{a(t+a )}(\comega)\subset d^aC^{t }(\comega)\text{ for
}t=\min\{2a,1,\tau -a\}-\varepsilon ;\text{ small }\varepsilon >0.
\end{equation}

$4^\circ$ In particular, the cokernels $N'_{\overline\lambda }$ of
$P_{D,2}-\lambda $, $\lambda \in\Sigma $, have the regularity stated
in $3^\circ$, and both kernels and cokernels are in all cases contained in $\dot C^a(\comega)$.
\end{thm}

It is understood here that $\varepsilon $ is chosen such that $\tau -a$, $t$ and
$t-a$ are not integer.

It is perhaps interesting to recall from \cite{G23} that the
identification of eigenspaces for different $q$ was
easy to show for $q\ge 2$, but demanded some effort for $q<2$, and vice
versa for the analysis of cokernels. Moreover, extra efforts were
necessary for results on the adjoint in the nonsymmetric case (when
$\tau <\infty $) because of
remainder terms in the pseudodifferential calculus. More details in \cite{G23}.

When $\tau \ge 3a$, the result in \eqref{eq:3.10} is only  an
$\varepsilon $ weaker than the result obtained for the $C^\infty $-case in
\cite{G15b} Th.\ 2.3.

\section{Weyl asymptotics for nonsmooth operators over smooth domains}\label{sec4}

It was shown in \cite{G15b} in the smooth case ($\tau =\infty $) that the eigenvalues of
$P_{D,2}$ (when selfadjoint), and more generally
the $s$-numbers, satisfy a Weyl asymptotic estimate where the
coefficient is determined by an integral formed of the principal symbol over
the domain. This will now be generalized to nonsmooth cases. In the
present paper we restrict the attention to eigenvalues in
symmetric cases (in the Hilbert space setting), where
there are some convenient auxiliary tools.

When $P$ is symmetric, the realization $ P_{D,2}$ is
 selfadjoint in $L_2(\Omega )$, and the sectorial set $M$ in \eqref{eq:3.6} is simply $M=[c_0-\beta ,\infty [\,$. Henceforth we drop the subscript 2,
 writing $P_{D,2}=P_D$. Here $P_{D}+\beta $ is
 a positive, bijective operator from $D(P_D)$ to $L_2(\Omega )$.

As already noted, the spectrum $\Sigma $ of $P_D$ is discrete. 
The
spectrum of $P_{D}+\beta $ is a sequence of positive eigenvalues 
(repeated according to multiplicities)
\begin{equation}\label{eq:4.1}
0<\lambda _1\le \lambda _2\le \dots \to \infty .
\end{equation}
The inverses of the $\lambda _j$ are the nonzero eigenvalues $\mu _j$
of the compact inverse of $P_{D}+\beta $. In the Appendix, Lemma
\ref{Lemma A.1},
we recall the equivalent formulations of asymptotic formulas for a
compact selfadjoint injective nonnegative operator  and its inverse.

For the extension of the Weyl formula from the smooth case, we shall
use some
well-known properties of $s$-numbers (also called singular values):

When $B$ is a compact linear operator from a Hilbert space $H$ to
another $H_1$, the $s$-numbers $s_j(B)$ (singular values) are defined as the numbers  $ s_j(B)=\mu
_j(B^*B)^\frac12$, where $\mu _j(B^*B)$ denotes the $j$-th positive
eigenvalue of $B^*B$, arranged nonincreasingly and repeated according
to multiplicities. A standard reference is the book of Gohberg and
Krein \cite{GK69}.
For $p>0$, the so-called weak Schatten class $\mathfrak
S_{p,\infty }(H,H_1)$ (with notation as in \cite{G14}),
consists of the compact operators $B$ such that
\begin{equation}\label{eq:4.2}
s_j(B)\le Cj^{-1/p}\text{ for all }j
.
\end{equation}
(The indication $(H,H_1)$ is replaced by $(H)$ if $H=H_1$; it can be
omitted when it is clear from the context.) The convention for $p$
stems from the fact that $\mathfrak
S_{p,\infty }$ is close to the standard Schatten class ${\mathfrak S}_{p } $
consisting of the operators $B$ with $(s_j(B))_{j\in{\mathbb N}}\in
\ell_p({\mathbb N})$. In fact, $\mathfrak
S_{p,\infty }\subset {\mathfrak S}_{p+\varepsilon } $ for any $\varepsilon
>0$. They are linear spaces. We recall some useful properties:
\begin{align}\label{eq:4.3}
&\mathfrak S_{p,\infty }\cdot \mathfrak S_{q,\infty }\subset \mathfrak S_{r,\infty
                             }, \text{ where }r^{-1}=p^{-1}+q^{-1},\\ \nonumber
&s_j(B^*)=s_j(B),\quad
s_j(EBF)\le \|E\|s_j(B)\|F\|,
\end{align}
when $E\colon H_1\to H_3$ and $F\colon H_2\to H$ are bounded linear maps between Hilbert
spaces. Moreover, there are perturbation rules:

\begin{lem}\label{Lemma 4.1}
$1^\circ$ If $s_j(B)j^{1/p}\to C_0$ and $s_j(B')j^{1/p}\to 0$ for
$j\to\infty$, then $s_j(B+B')j^{1/p}\to C_0$ for $j\to\infty$.

$2^\circ$ If $B=B_M+B'_M$ for each $M\in\mathbb N$, where
$s_j(B_M)j^{1/p}\to C_M$ for $j\to\infty$ and $s_j(B'_M)j^{1/p}\le
c_M$ for
$j\in\mathbb N$, with $C_M\to C_0$ and $c_M\to 0$ for
$M\to\infty$, then $s_j(B)j^{1/p}\to C_0$ for $j\to\infty$.

\end{lem}

The statement in $1^\circ$ is the Weyl-Ky Fan theorem (cf.\ e.g.\
\cite{GK69} Th.\ II 2.3), and $2^\circ$ is a refinement shown in
\cite{G84} Lemma 6.2.$2^\circ$. More details in \cite{G96} Sect.\
A.6. One also defines singular values for unbounded operators
$A$ with discrete spectrum
(including inverses of operators $B$ as above); we shall use the
notation $r_j(A)=\lambda _j(A^*A)^\frac12$ as in \cite{G96}. By an abuse
of notation, they were called $s_j$ again in \cite{G15b}.
\medskip

In the case where the symbol of $P$ depends smoothly on $x$, we have
from \cite{G15b} Th.\ 2.7:

\begin{thm}\label{Theorem 4.2} \cite{G15b} Assume Hypothesis \ref{Hypothesis 3.1} with $\tau
=\infty $, and let $\Omega $ be  bounded and $C^{\infty  }$. The singular values $r_j(P_D)$ of $P_D$ (the eigenvalues of $(P^*_D
P_{D})^\frac12$), as well as the singular values of $P_D+b$ for any $b\in\C$, have the asymptotic behavior
\begin{equation}\label{eq:4.4}
r_j(P_{D}+b)=C(P,\Omega )j^{2a/n}+o(j^{2a/n})
\text{ for }j\to\infty ,
\end{equation}
where $C(P,\Omega )$ is
defined from the principal symbol $p_{0}(x,\xi )$ by: 
\begin{equation}\label{eq:4.5}
C(P,\Omega )=C'(P,\Omega )^{-2a/n},\quad C'(P,\Omega )=\tfrac 1{n(2\pi
  )^n}\int_{\Omega }\int_{|\xi |=1}|p_{0}(x,\xi )|^{-n/2a}\,d\omega
(\xi )dx.
\end{equation} 
\end{thm}

We remark here that the proof of Theorem 2.7 in \cite{G15b} first shows that the asymptotic
behavior holds for $P_D+b$ when $P_D+b$ is invertible, and then
concludes it for $P_D$ by a perturbation argument. (One can here take
$b=\beta $ in \eqref{eq:3.3}.) An analogous perturbation argument allows to include
arbitrary $b$.

Now consider nonsmooth symbols. When $\tau <\infty $, we approximate the symbol $p(x,\xi )$   in $C^{\tau '}S^{2a}$ ($\tau '<\tau $ close to $\tau $) by smooth
symbols $p_k(x,\xi )$, $k\to\infty $, obtained by
convolutions in $x$ with a resolution of the identity as in
\cite{G14};
the $p_k$ are likewise strongly elliptic and even. They give rise to operators
$  P_k$ on $\rn$ and Dirichlet realizations $P_{k,D}$ on $\Omega $. By \cite{G14} (2.18), the
operator norm of $ P-P_k$ in  $\mathcal  L( H^a(\rn),
H^{\,-a}(\rn))$ 
goes to 0, so we get from \eqref{eq:3.3} that
\begin{equation}\label{eq:4.6}
\operatorname{Re} (P_ku,u)\ge c/2\|u\|^2_{ H^a(\rn)}-(\beta +1)\|u\|^2_{L_2(\rn )} 
\end{equation}
for sufficiently large $k$; discarding the first terms, we can assume
that \eqref{eq:4.6} holds for all $k$.

For $b\in\R$, denote $P_D+b=P_{D,b}$ and $P_{k,D,b}=P_{k,D}+b$; then
these operators have positive lower bound and are invertible when $b\ge
\beta +1$. 

The asymptotic estimate of the eigenvalues of
$P_{D,b}^{-1}$ will be obtained by approximation from the corresponding
result for the $s$-numbers of  $P_{k,D,b}^{-1}$, that holds according to
Theorem \ref{Theorem 4.2}.
For that purpose, we need to estimate the $s$-numbers of
\begin{equation}\label{eq:4.7}
P_{D,b}^{-1}-P_{k,D,b}^{-1}=P_{k,D,b}^{-1}(P_{k,D,b}-P_{D,b})P_{D,b}^{-1}.
\end{equation}

Since $H^{a(2a)}(\comega)$ is only for $0<a<\frac12$ known to be
contained in a $2a$-order Sobolev space, the embedding of
$H^{a(2a)}(\comega)$ into $L_2(\Omega )$ needs special considerations:

\begin{lem}\label{Lemma 4.3} When $\Omega \subset \rn$ is smooth and bounded,
and $T$ is a linear operator in $L_2(\Omega )$ that is continuous from $L_2(\Omega )$ to
$H^{a(2a)}(\comega )$, then its $s$-numbers satisfy
\begin{equation}\label{eq:4.8}
s_j(T)\le C_1  \|T\|_{\mathcal L(L_2, H^{a(2a)})} j^{-2a/n}, \text{ all }j, 
\end{equation}
where $C_1$ is a constant depending only on $\Omega $ and $a$.
\end{lem}
 
\begin{proof} This follows since we know from \cite{G15b} that for $P_{\Delta } =(-\Delta )^a$,
the inverse of the Dirichlet realization is a homeomorphism $P_{\Delta ,
D}^{-1}$  from
$L_2(\Omega )$ to $H^{a(2a)}(\comega )$ satisfying an estimate
(as a consequence of Theorem \ref{Theorem 4.2}):
\begin{equation*}
s_j(P_{\Delta,D}^{-1})\le  C_0j^{-2a/n},
\text{ all }j.
\end{equation*}
Then, since $T$ is
the composition of $P_{\Delta
,D}^{-1}$ and the bounded operator $B=P_{\Delta
,D}T$  in $L_2(\Omega )$,
the property follows for $T$ by the last rule in \eqref{eq:4.3}:
\begin{align*}
s_j(T)&=s_j( P_{\Delta
,D}^{-1}B)\le C_0j^{-2a/n}\|B\|\\
&\le
C_0\|P_{\Delta
,D}\|_{\mathcal L(H^{a(2a)},L_2)}\|T\|_{\mathcal L(L_2,
H^{a(2a)})}j^{-2a/n}= C_1\|T\|_{\mathcal L(L_2,
H^{a(2a)})}j^{-2a/n}. 
\end{align*}

\end{proof}

\begin{pro}\label{Proposition 4.4} Let $b\ge \beta +1$.
The operator norm of the  difference between
$P_{D,b}^{-1}$ and $P_{k,D,b}^{-1}$ satisfies
\begin{equation}\label{eq:4.9}
\|P_{D,b}^{-1}-P_{k,D,b}^{-1}\|_{\mathcal
L(L_2, H^{a(2a)})}\to 0 \text{ for }k\to\infty ,
\end{equation}
and its $s$-numbers
satisfy inequalities
\begin{equation}\label{eq:4.10}
s_j(P_{D,b}^{-1}-P_{k,D,b}^{-1})/j^{-2a/n}\le
C_k\text{ for all }j,
\text{ where }C_k\to 0 \text{ for }k\to \infty .
\end{equation}
\end{pro}

\begin{proof} We use the representation \eqref{eq:4.7}.  Since
$D(P_{D,b})=H^{a(2a)}(\comega)$,
\begin{equation}\label{eq:4.11}
\|P_{D,b}^{-1}\|_{\mathcal
L(L_2, H^{a(2a)})}\le C'.
\end{equation}

To the difference
 $P_{k,D,b}-P_{D,b}$ we apply
 Theorem 5.11 of \cite{AG23}, supplied with norm considerations as in
 Theorem 5.13 there. Namely, the boundedness of
 $P_{k,D,b}-P_{D,b}$ follows from
 Theorem 5.11, with $\mu =a$, $m=2a$, $s=0$, $q=2$, and $\tau $
 replaced by $\tau '< \tau $ close to $\tau $; and for the pieces localized as in the proof of
 Theorem 5.11, the convergence
 to zero in operator norm follows from the convergence to zero of a
 sufficiently high-numbered symbol seminorm, as shown
 in Theorem 5.13 there. Thus 
\begin{equation}\label{eq:4.12}
\|P_{D,b}-P_{k,D,b}\|_{\mathcal
L( H^{a(2a)},L_2)}\to 0 \text{ for }k\to\infty .
\end{equation}

It follows in particular that
\begin{equation*}
P_{k,D,b}P_{D,b}^{-1}\to I \text{ in }\mathcal
L(L_2(\Omega )) \text{ for }k\to\infty ,
\end{equation*}
so by a Neumann series argument, the inverse
$P_{D,b}P_{k,D,b}^{-1}$ likewise has 
operator norm in $L_2(\Omega )$
bounded with respect to  $k$, and hence
\begin{equation}\label{eq:4.13}
\|P_{k,D,b}^{-1}\|_{\mathcal
L(L_2, H^{a(2a)})}\le C''\text{ for all }k.
\end{equation}
Applying \eqref{eq:4.11}, \eqref{eq:4.12} and \eqref{eq:4.13} to the factors in \eqref{eq:4.7}, we
conclude \eqref{eq:4.9}, and then \eqref{eq:4.10} follows by Lemma
\ref{Lemma 4.3}.

\end{proof}

\begin{thm}\label{Theorem 4.5} Assume Hypothesis \ref{Hypothesis 3.2}, let $\Omega
\subset \rn$ be  bounded and $C^\infty $, and let $P_{D}$ be the Dirichlet realization
of $P$ in $L_2(\Omega )$. Let $b\in\R$ be such that $P_{D,b}=P_D+b$
has positive lower bound. Then the
eigenvalues $\mu _j(P_{D,b}^{-1})$ of its inverse satisfy:
\begin{equation}\label{eq:4.14}
\mu
_j(P_{D,b}^{-1})=C(P,\Omega )^{-1}j^{-2a/n}+o(j^{-2a/n})\text{
for }j\to \infty ,
\end{equation}
where $C(P,\Omega )$ is defined in \eqref{eq:4.5}.
Equivalently, the eigenvalues $\lambda  _j(P_{D})$
satisfy:
\begin{equation}\label{eq:4.15}
\lambda 
_j(P_{D})=C(P,\Omega )j^{2a/n}+o(j^{2a/n})\text{
for }j\to \infty .
\end{equation}

\end{thm}

\begin{proof} Define the approximating operators
$P_{k,D,b}$ as described above. We have from Theorem 4.2 that
their $s$-numbers satisfy:
\begin{equation}\label{eq:4.16}
s_j(P_{k,D,b}^{-1})=C(P_k,\Omega )^{-1}j^{-2a/n}+o(j^{-2a/n})
\text{ for }j\to \infty ,  
\end{equation}
where the constants are defined as in \eqref{eq:4.5}. For the differences
$P_{D,b}^{-1}-P_{k,D,b}^{-1}$ we have
the estimates of $s$-numbers in \eqref{eq:4.10}. Clearly,
$C(P_k,\Omega )\to C(P,\Omega )$
for $k\to \infty $ (by the convergence of the principal symbols). Then the result follows by use of Lemma 4.1 $2^\circ$,
noting that $\mu
_j(P_{D,b}^{-1})=s_j(P_{D,b}^{-1})$.

The last statement follows in view of Lemma A.1 applied to
$P_{D,b}$. Clearly, the asymptotic estimate \eqref{eq:4.15} is equivalent with
the estimate where $P$ is replaced by $P+b$.
\end{proof}

The result extends to various perturbations of $P$. One is to add to
$P$ an operator $P'$ with lower positive order and suitable mapping
properties. Another is to add to $P$ an $L_2$-bounded operator $P''$.

\begin{thm}\label{Theorem 4.6} Assume Hypothesis \ref{Hypothesis 3.2} and let $\Omega
\subset \rn$ be  bounded and $C^{\infty  } $.
Let $P'$ and $P''$ be  symmetric
operators that  map continuously, for some $0<\delta \le a$,
\begin{align}\label{eq:4.17}
P'\colon  H^{t }(\rn)&\to H^{t-2a+\delta }(\rn)\text{ for }a-\delta \le
t\le 2a ,
\text{ if } a<\tfrac12;\\ \label{eq:4.18}
P'\colon  H^{t }(\rn)&\to  H^{t-a-\frac12+\delta }
(\rn )\text{ for }\tfrac12-\delta \le t\le a+\tfrac12,
\text{ if } a\ge \tfrac12;\\ \label{eq:4.19}
  P''\colon  L_2(\rn)&\to  L_2(\rn);
                       \end{align}
$P''$ can for example be the multiplication by a function $V(x)\in
L_\infty (\rn,\R)$.

$1^\circ$ The Dirichlet realization $\tilde P_D$ of $\tilde
P=P+P'+P''$ defined from the sesquilinear form $\tilde s$ on $\dot
H^a(\comega)$ extending
\begin{equation}\label{eq:4.20}
\tilde s(u,v)=\int_\Omega (P+P'+P'')u\, \bar v\, dx \text{ on }C_0^\infty (\Omega ),
\end{equation}
is selfajoint lower bounded and has $D(\tilde
P_D)=H^{a(2a)}(\comega)$.

$2^\circ$  
The eigenvalues of $\tilde
P_D$ satisfy:
\begin{equation}\label{eq:4.21}
\lambda _j(\tilde P_D)=C(P,\Omega )j^{2a/n}+o(j^{2a/n})\text{
for }j\to \infty ,
\end{equation}
with the constant $C(P,\Omega )$ defined from $p_0$ as in  \eqref{eq:4.5}.
\end{thm}

\begin{proof} $1^\circ$.
We first show how to handle the addition of $P'$. Let $s'(u,v)$ be the
form
$s'(u,v)=(P'u,v)$ for $u,v\in C_0^\infty (\Omega )$.

\emph{The case $a<\frac12$.}  By use of the order-reducing operators
$\Lambda _+^{(t)}$ and their adjoints $\Lambda _{-,+}^{(t)}\equiv r^+\Lambda _-^{(t)}e^+$ recalled in \eqref{eq:2.7}, we have that for $u,v\in
C_0^\infty (\Omega )$, 
\begin{align}\label{eq:4.22}
|(P'u,v)_{L_2(\Omega )}|&=|(r^+P'u,\Lambda _+^{(-a)}\Lambda
 _+^{(a)}v)_{L_2(\Omega )}|=|(\Lambda _{-,+}^{(-a)}r^+P'u,\Lambda
 _+^{(a)}v)_{L_2(\Omega )}|\\ \nonumber
 &\le c\|\Lambda _{-,+}^{(-a)}r^+P'u\|_{L_2(\Omega )}\|\Lambda
 _+^{(a)}v\|_{L_2(\Omega )}|\le c'\|u\|_{\dot H^{a-\delta
 }(\comega)}\|v\|_{\dot H^a(\comega)};
 \end{align}
where we used \eqref{eq:4.17} for $t=a-\delta $ (here $\Lambda _{-,+}^{(-a)}$
lifts $\ol H^{\,-a}(\Omega )$ to $L_2(\Omega )$). So the form
$s'(u,v)$ extends to a form $s'(u,v)$
continuous on $\dot H^a(\comega)$, and moreover, by interpolating
$\dot H^{a-\delta }(\comega)$ between $\dot H^a(\comega)$ and $L_2(\Omega )$ we can deduce that for any $\varepsilon >0$,
\begin{equation}\label{eq:4.23}
|s'(u,u)|\le \varepsilon \|u\|^2_{\dot H^a(\comega)}+C(\varepsilon
)\|u\|^2_{L_2(\Omega )}.
\end{equation}
Taking $\varepsilon \le c_0/2$, we conclude that $s(u,v)+s'(u,v)$
satisfies a coerciveness estimate like \eqref{eq:3.3} over $\Omega $:
\begin{equation*}
s(u,u)+s'(u,u)\ge c_0/2\|u\|^2_{\dot
H^a(\comega)}-\beta '\|u\|^2_{L_2(\Omega )}, \text{ for }u\in \dot H^a (\comega),
\end{equation*}
 with a  larger $\beta '$. Then the
variational construction carries through to define the Dirichlet
realization $(P+P')_D$, selfadjoint and with lower bound $\ge
c_0/2-\beta '$.

We shall show that its domain (defined by \eqref{eq:3.5} for $P+P'$) is again $H^{a(2a)}(\comega)$, in the
present case equal to $\dot H^{2a}(\comega)$: When $u\in \dot
H^{2a}(\comega)$, then we know that $r^+Pu\in L_2(\Omega )$, and we
have  $r^+P'u\in L_2(\Omega )$ by \eqref{eq:4.17} for $t=2a$. Conversely, let
$u\in \dot H^a(\comega )$ satisfy
\begin{equation}\label{eq:4.24}
r^+Pu+r^+P'u=f, \quad f\in L_2(\Omega ) .
\end{equation}
By \eqref{eq:4.17} with $t=a$, $f-r^+P'u\in \ol H^{\min\{0,-a+\delta \} }(\Omega )$, so by
the regularity of the Dirichlet problem for $P$, $u\in \dot
H^{\min\{2a,a+\delta \} }(\comega)$, cf.\ Theorem \ref{Theorem 3.3}. If $\delta = a$, we are through. If not,
we use \eqref{eq:4.17} again, now with $t=a+\delta $, to see that $f-r^+P'u\in
\ol H^{\min\{0,-a+2\delta \} }(\Omega )$, hence $u\in \dot
H^{\min\{2a,a+2\delta \} }(\comega)$. In finitely many steps of this
kind, the conclusion $u\in \dot H^{2a}(\comega)$ is reached.

\emph{The case $a\ge \frac12$.} Here we find \eqref{eq:4.22} by use of \eqref{eq:4.18}
with $t=\frac12-\delta $. Then \eqref{eq:4.23} can be concluded, leading to a
selfadjoint lower bounded Dirichlet realization $(P+P')_D$.

To show that the domain is  $H^{a(2a)}(\comega)$, we proceed as
follows: Since $H^{a(2a)}(\comega)\subset \dot H^{a+\frac12-\delta
}(\comega)$, $r^+(P+P')$ maps $H^{a(2a)}(\comega)$ into $L_2(\Omega )$
by \eqref{eq:4.18} for $t=a+\frac12$. Conversely, consider \eqref{eq:4.24}.
By \eqref{eq:4.18} with $t=a$, $f-r^+P'u\in \ol H^{\min\{0,-\frac12+\delta \}
}(\Omega )$, so by Theorem \ref{Theorem 3.3}, $u\in 
H^{a(\min\{2a,2a-\frac12+\delta \}) }(\comega)$. If $\delta \ge
\frac12$, we are through. If not, we observe that
$H^{a(\min\{2a,2a-\frac12+\delta \}) }(\comega)\subset H^{a(a+\delta
\} }(\comega)$ since $a\ge \frac12$; here
\begin{equation*}
 H^{a(a+\delta
) }(\comega)=\Lambda _+^{(-a)}e^+\ol H^\delta (\Omega )=\Lambda
_+^{(-a)} \dot H^\delta (\comega )=\dot H^{a+\delta }(\comega),
\end{equation*}
since $\delta <\frac12$. By \eqref{eq:4.18} with $t=a+\delta $, $r^+P'u\in \ol
H^{\,-\frac12+2\delta }(\Omega )$, hence 
 $f-r^+P'u\in
\ol H^{\min\{0,-\frac12+2\delta \} }(\Omega )$. It follows that $u\in 
H^{a(\min\{2a,2a-\frac12+2\delta \}) }(\comega)$. In finitely many steps of this
kind, we reach the conclusion $u\in  H^{a(2a)}(\comega)$.

Adding $P''$ to $P+P'$ is much easier, since $r^+P''e^+$ defines a 
bounded map in $L_2(\Omega )$. When we add $(P''u,v)$ to
$s(u,v)+s'(u,v)$, we can easily verify the $\dot H^a(\Omega
)$-continuity and coerciveness for the resulting form $\tilde s(u,v)$. Since
$H^{a(2a)}(\comega)\subset L_2(\Omega )$, $\tilde P$ is continuous
from $H^{a(2a)}(\comega)$ to $ L_2(\Omega )$. In the equation
\begin{equation*}
r^+(P+P'+P'')u=f,\quad f\in L_2(\Omega )
,\; u\in \dot H^a(\comega),
\end{equation*}
the right-hand side remains in $L_2(\Omega )$, when $r^+P''u\in L_2(\Omega
)$ is moved there. So we can conclude that $D(\tilde
P_{D})=H^{a(2a)}(\comega)$.

$2^\circ$. 
We use Theorem
\ref{Theorem 4.2} to derive the Weyl asymptotics formula for $\tilde P_D$ from
 that of $P_D$ by perturbation:

Take $b$ so large that both $P_D+b$ and $\tilde P_D+b$ have positive
lower bound. Then
\begin{equation*}
(P_D+b)^{-1} - (\tilde P_D+b)^{-1}=(\tilde P_D+b)^{-1}(r^+P'+r^+P'')( P_D+b)^{-1}.
\end{equation*}

Here, by the assumptions on $P'$, $r^+P'$ maps $H^{a(2a)}(\comega)$
into $\ol H^\delta (\Omega )$ if $a<\frac12$ and into e.g.\ $\ol
H^{\delta/2}(\Omega )$ if $a\ge \frac12$, so $r^+P'( P_D+b)^{-1}$ maps
$L_2(\Omega )$ to $\ol H^{\delta/2}(\Omega )$ and hence belongs to the
weak Schatten class ${\mathfrak S}_{n/(\delta /2),\infty  } $. As for the
contribution from $P''$, we use its continuity in $L_2(\Omega )$ and the
last rule in \eqref{eq:4.3} to see that $r^+P''( P_D+b)^{-1}$ is in  ${\mathfrak
S}_{n/(2a),\infty  } $. Here $\delta /2< 2a$. Then by the 
product rule in  \eqref{eq:4.3},
\begin{align*}
(\tilde P_D+b)^{-1}(r^+P'+r^+P'')( P_D+b)^{-1}&\in {\mathfrak S}_{n/(2a),\infty
} {\mathfrak S}_{n/(\delta /2),\infty  }\\& \subset {\mathfrak S}_{n/(2a+\delta /2),\infty  }. 
\end{align*}
Finally, the first rule in Lemma \ref{Lemma 4.1} implies the asserted asymptotics for
$\tilde P_D$. 

\end{proof}

$1^\circ$ also holds when $\Omega $ is $C^{1+\tau }$; here the use of
the order-reducing operators $\Lambda _+^t$ to show \eqref{eq:4.23} is
replaced by the use of $\Xi _+^t=\Op((\ang{\xi '}+i\xi _n)^t)$
 in local coordinates.

The hypotheses on $P'$ are satisfied if $P'$ is a pseudodifferential
operator with symbol in $C^\tau S^{a+\min\{a,\frac12\}-\delta
}_{1,0}(\rn\times\rn)$ (recall $0<\delta \le a$). In particular, $P'$
can be of order $a$ with symbol in $C^\tau S^{a}_{1,0}(\rn\times\rn)$. And, as mentioned, $P''$ can be a bounded real potential $V(x)$.

Note that only the action of $r^+(P'+P'')$ on functions supported in
$\comega$ really enters, so the requirements \eqref{eq:4.17}--\eqref{eq:4.19} could be
replaced by properties referring only to $\comega$.

\section{Nonsmooth operators over nonsmooth sets}\label{sec5}

The asymptotic formula \eqref{eq:4.21} can be extended to nonsmooth domains by a simple argument
in the selfadjoint case. As a generalization of the notion of ``contented'' in Reed and Simon
\cite{RS78} p.\ 271, we define:

\begin{defn}\label{Definition 5.1} We shall say that a bounded open set $\Omega \subset\rn $ is $C^\infty $-contented, when
there exist two sequences of $C^\infty $-domains
$\Omega _{\operatorname{in},j}\subset \Omega , j\in{\mathbb N}$, and  
$\Omega _{\operatorname{out},k}\supset \Omega, k\in{\mathbb N} $, such that
\begin{equation}\label{eq:5.1}
\lim_{j\to \infty }\operatorname{vol}(\Omega
_{\operatorname{in},j})=\lim_{k\to \infty }\operatorname{vol}(\Omega
_{\operatorname{out},k})=\operatorname{vol}(\Omega ).
\end{equation}
\end{defn}
 
This holds for example for a Lipschitz domain (i.e.,
$C^{0,1}$-domain):

\begin{lem}\label{Lemma 5.2} When $\Omega $ is bounded with a
$C^{0,1}$-boundary, then $\Omega $ is $C^\infty $-contented.
\end{lem}

\begin{proof} This can be shown by an elaboration of the explanation
given in Daners \cite{D08} Prop.\ 8.2.1. For the approximation from
inside, let
\begin{equation}\label{eq:5.2}
  V_j=\{x\in\Omega \mid \operatorname{dist}(x,\partial\Omega )> 1/j\},
\end{equation}
then $\{V_j\}_{j\in{\mathbb N}}$ is a nested sequence of open sets such
that $\overline V_j\subset V_{j+1}\subset \overline V_{j+1}\subset
\Omega $ for all $j$. Since $\Omega $ is $C^{0,1}$,
$\operatorname{vol}(\Omega \setminus V_j)\to 0$ for $j\to \infty
$. Namely, every boundary point $x$ has a bounded neighborhood
$W_{x}$ such that, after a translation and rotation, 
 $W_{x}\cap \Omega =W_{x}\cap \rn_\zeta$ ($\rn_\zeta =\{x=(x',x_n)\mid x_n>\zeta (x')\}$) for a
$C^{0,1}$-function $\zeta (x')$. Here $(W_{x}\cap(\Omega \setminus V_j))\cap
\rn_\zeta \subset W_{x} \cap(\rn_{\zeta +C/j}\setminus \rn_\zeta )$
for a constant $C$, so that the volume goes to $0$ for $j\to\infty
$. Since the boundary can be covered by finitely many such sets
$W_{x}$, $\operatorname{vol}(\Omega \setminus V_j)\to 0$ for
$j\to\infty $.

Now as shown in \cite{D08} Prop.\ 8.2.1, we can for each $j$ find a
$C^\infty $-set $U_{\operatorname{in},j}$ with $\overline V_j\subset
U_{\operatorname{in},j}\subset \overline
U_{\operatorname{in},j}\subset V_{j+1}$ by choosing a cutoff function
$\psi _j\in C_0^\infty (V_{j+1},[0,1])$ that equals $1$ on $V_j$; then $\psi
_j$ has by Sard's lemma (cf.\ e.g.\ Hirsch \cite{H76} Th.\ 3.1.3) a regular value $t_j\in \,]0,1[\,$ such that
$U_{\operatorname{in},j}=\{x\mid \psi (x)>t_j\}$ is $C^\infty $.

This defines the sets $U_{\operatorname{in},j}$ approximating $\Omega
$ from inside, and a similar study of $(\rn\setminus \comega)\cap B_R$,
for a large $R$ such that the ball $B_B=\{|x|<R\}$ contains $\comega$, gives
an approximating  sequence $U_{\operatorname{out},k}$ from outside.

\end{proof}

Let $\tilde P$ be as in Theorem \ref{Theorem 4.6}. Recall that the eigenvalues of
$\tilde P_{D}$ can be
described as the Rayleigh quotients, where $H=L_2(\Omega )$,
$V=\dot H^a(\comega)$, $X$ denotes a linear space:
\begin{equation}\label{eq:5.3}
\lambda _j(\tilde P_{D})=\underset{X\subset H,\dim X< j}\max\; \underset{v\in V,v\perp
X}\min\,\frac{\tilde s(v,v)}{\|v\|^2_H},
\end{equation}
and $\tilde s(u,v)$ is the sesquilinear form on $\dot H^a(\comega)$
defined by \eqref{eq:4.20}.

There are likewise defined Dirichlet realizations $\tilde P_{D,\Omega _{\operatorname{in},l}}$ of $\tilde P$ on the sets
$\Omega _{\operatorname{in},l}$;
with associated forms
\begin{equation}\label{eq:5.4}
\tilde s_{\Omega
_{\operatorname{in},l}}(v,v)=\int_{\Omega
_{\operatorname{in},l}}\tilde Pv\,\bar v\,dx\text{ on }\dot H^a(\comega
_{\operatorname{in},l}),
\end{equation}
 coinciding with $\tilde s(v,v)$ for $v\in \dot H^a(\comega
 _{\operatorname{in},l})$. Dirichlet realizations
 $\tilde P_{D,\Omega _{\operatorname{out},l}}$ on $\Omega
 _{\operatorname{out},l}$ are similarly defined.
 Based on the description of the eigenvalues as Rayleigh quotients
 \eqref{eq:5.3}, we can compare the eigenvalues of
 $\tilde P_{D,\Omega _{\operatorname{in},l}}$ and
 $\tilde P_{D}$, showing that the eigenvalues of
$\tilde P_{D,\Omega _{\operatorname{in},l}}$ are larger
than those of $\tilde P_{D}$ for each $j$.
This follows by application of a general well-known comparison
 property Proposition \ref{Proposition A.2}, that we include in the Appendix (with a
 proof for the convenience of the reader).

We can then show:

\begin{thm}\label{Theorem 5.3} Let $\tilde P=P+P'+P''$ be as in
  Theorem \ref{Theorem 4.6}.
Let $\Omega \subset \rn$ be bounded and $C^\infty $-contented; e.g.\
 a Lipschitz domain.
Then the eigenvalues of the Dirichlet realization $\tilde P_{D}$
of $\tilde P$ in $L_2(\Omega )$ satisfy:
\begin{equation}\label{eq:5.5}
\lambda 
_j(\tilde P_{D})=C(P,\Omega )j^{2a/n}+o(j^{2a/n})\text{
  for }j\to \infty ,
\end{equation}
where $C(P,\Omega )$ is defined by \eqref{eq:4.5}.

\end{thm}

\begin{proof} The variational construction of $\tilde P_D$ also works when
$\Omega $ is $C^{0,1 }$ (and even for less smooth domains), giving a
selfajoint lower bounded operator with compact resolvent. We apply the comparison principle recalled in
Proposition \ref{Proposition A.2} with $H_1=L_2(\Omega _{\operatorname{in},l})$,
$V_1=\dot H^a(\overline\Omega _{\operatorname{in},l})$,
$s_1(u,v)=\tilde s_{\Omega _{\operatorname{in},l}}(u,v)$ (cf.\ \eqref{eq:5.4}), and  $H_2=L_2(\Omega)$,
$V_2=\dot H^a(\comega)$, $s_2(u,v)=\tilde s(u,v)$ (cf.\ \eqref{eq:4.20}); this gives that 
\begin{equation}\label{eq:5.6}
\lambda _j(\tilde P_{D,\Omega _{\operatorname{in},l}})\ge
\lambda _j(\tilde P_{D}), \text{ all }j,l.
\end{equation}
A similar comparison between $H_1=L_2(\Omega)$,
$V_1=\dot H^a(\comega)$, $s_1(u,v)=\tilde s(u,v)$, and $H_2=L_2(\Omega _{\operatorname{out},l})$,
$V_2=\dot H^a(\overline\Omega _{\operatorname{out},l})$,
$s_2(u,v)=\tilde s_{\Omega _{\operatorname{out},l}}(u,v)$, gives
\begin{equation}\label{eq:5.7}
\lambda _j(\tilde P_{D})\ge
\lambda _j(\tilde P_{D,\Omega _{\operatorname{out},l}}),  \text{ all }j,l.
\end{equation}
Thus we have for each $l$, applying Theorem \ref{Theorem 4.6} to $\tilde P_{D,\Omega _{\operatorname{in},l}}$,
\begin{equation}\label{eq:5.8}
\lim\sup_{j\to\infty }\lambda _j(\tilde P_{D})j^{-2a/n}\le
\lim\sup_{j\to\infty }\lambda _j(\tilde P_{D,\Omega
  _{\operatorname{in},l}})j^{-2a/n}=C(P,\Omega_{\operatorname{in},l}).
\end{equation}
By \eqref{eq:4.5} and \eqref{eq:5.1},  
\begin{equation*}
|C(P,\Omega
) -C(P,\Omega
_{\operatorname{in},l})|\le C'\operatorname{vol}(\Omega \setminus
\Omega _{\operatorname{in},l})\to 0 \text{ for }l\to \infty . 
\end{equation*}
Then letting $l\to\infty $ in \eqref{eq:5.8}, we find that
\begin{equation*}
\lim\sup_{j\to\infty }{\lambda _j(\tilde P_{D})}j^{-2a/n}
\le C(P,\Omega
) . 
\end{equation*}
It is shown similarly by use of \eqref{eq:5.7} 
that
\begin{equation*}
\lim\inf_{j\to\infty }{\lambda _j(\tilde P_{D})}j^{-2a/n}\ge C(P,\Omega
) , 
\end{equation*}
and \eqref{eq:5.5} follows.
\end{proof}

\begin{exam}\label{Example 5.4} Theorem \ref{Theorem 5.3} applies e.g.\ to $(-\Delta
)^a+(-\Delta )^{a'}+V(x)$, with  $0<2a'<
\min\{2a,a+\frac12\}$ and $V\in L_\infty(\rn,\R) $; here  $P=(-\Delta )^a$, $P'=(-\Delta
)^{a'}$ and $P''=V$. For a more general
example, we can replace $-\Delta $ by  a
selfadjoint strongly elliptic nonnegative second-order differential operator
$L=\sum_{|\alpha |\le 2}a_\alpha  (x)\partial^\alpha $ with $C^\infty
$-coefficients, defining the fractional power  $L^a$ by Seeley's
construction \cite{S67}. In particular, $P$ can
be the fractional
power
$\bigl((i\nabla+A(x))^2+m^2\bigr)^a$ of a magnetic
Schr\"odinger operator 
$(i\nabla+A(x))^2+m^2$, $A\in C^\infty (\rn,\rn)$. For $P'$ we can take
a general
symmetric pseudodifferential operator of order $2a'$ with $C^\tau $-symbol.
\end{exam}

Besides the examples mentioned above, we  note that the theorem
applies with $P$ replaced by
$P_{\operatorname{Re}}=\frac12(P+P^*)$ for an operator $P$ satisfying Hypothesis
\ref{Hypothesis 3.1} with $\tau =\infty $.

\begin{rem}\label{Remark 5.5} Departing from Theorem \ref{Theorem 5.3}, one can now perturb $\tilde P_D$ by suitable
lower-order nonsymmetric operators to obtain similar estimates of
singular values  by use of the
principles recalled in the start of Section \ref{sec4}. As for taking $P$
itself nonsymmetric when $\Omega $ is nonsmooth, it seems that the
analysis of
singular values may
demand a more circumstantial effort. One possibility is to use coordinate changes from
$\Omega $ to neighboring smooth domains (in the spirit of \cite{G14}
Th.\ 6.2), based on the results on nonsmooth coordinate changes for
$\psi $do's in \cite{AG23}.
\end{rem}

\section {Appendix.  Rules for eigenvalues}\label{sec6}

For the study of eigenvalues we recall the elementary transition between the eigenvalue behavior of a
compact nonnegative injective operator $B$, the eigenvalue behavior of its
inverse $A$,
and the counting function for $A$ (taken up in detail e.g.\ in
\cite{G96} Lemma A.5):

\begin{lem}\label{Lemma A.1} Let $B$ be compact selfadjoint $\ge 0$ and
injective, with eigenvalues $\mu _j(B)$ going nonincreasingly to $0$,
let $A=B^{-1}$ with eigenvalues $\lambda _j(A)=\mu _j(B)^{-1}$ going nondecreasingly to $\infty $,
and let $N(t;A)$ denote the number of eigenvalues of $A$ in $[0,t]$
(all eigenvalues counted with  multiplicity). Let $C_0>0$, $p>0$. The
following three statements are equivalent:
\begin{align} \nonumber
\mu _j(B)&=C_0^{1/p}j^{-1/p}+o(j^{-1/p}) \text{ for }j\to\infty ,\\\label{eq:A.1}
  \lambda  _j(A)&=C_0^{-1/p}j^{1/p}+o(j^{1/p})\text{ for }j\to\infty ,\\
  \nonumber
N(t; A)&=C_0t^p +o(t^p) \text{ for }t\to\infty .
\end{align}

\end{lem}

The following is a well-known comparison principle for eigenvalues of
operators defined from sesquilinear forms (with notation as in e.g.\
\cite{G09} Sect.\ 12.4, after Lions-Magenes \cite{LM68}):

\begin{pro}\label{Proposition A.2} Let $(H_1,V_1,s_1)$ and $(H_2,V_2,s_2)$ be triples
giving rise to selfadjoint variational operators $A_1$ resp.\ $A_2$ by
the Lax-Milgram lemma. 
 Assume that $V_1\subset V_2$ with continuous injection, that
$H_1$ is a closed subspace of $H_2$, that the
injections of $V_i$ into $H_i$ are compact ($i=1,2$), and that $s_1(v,v)\ge
s_2(v,v)$ for $v\in V_1$. Then the eigenvalues $\lambda _j(A_1)$ and
$\lambda _j(A_2)$ of $A_1$ and $A_2$ (ordered nondecreasingly) satisfy
\begin{equation}\label{eq:A.2}
\lambda _j(A_1)\ge \lambda _j(A_2),\text{ for all }j\in\mathbb N.
\end{equation}
 \end{pro}

\begin{proof} The eigenvalues can be represented as Rayleigh quotients
  (with $v\ne 0$):
\begin{equation}\label{eq:A.3}
\lambda _j(A_i)=\underset{X\subset H_i,
\dim X< j}\max\;\underset{v\in V_i,v\perp
X}\min\,\frac{s_i(v,v)}{\|v\|^2_{H_i}},\quad i=1,2.
\end{equation}
Note that when $X_2$ is a finite dimensional subspace of $H_2$, and
$v\in H_1$, then $v\perp X_2\iff v\perp X_1$, where
$X_1=\Pi _{H_1}X_2$, orthogonal projection. All subspaces of $H_1$ of
dimension $\le j-1$ are obtained in the form $\Pi _{H_1}X$ when $X$ runs
through the subspaces of $H_2$ of dimension $\le j-1$.

For each $X\subset H_2$ of dimension $j-1$,
\begin{align} \nonumber
\min \{\,&s_1(v,v)/\|v\|_{H_1}^2\mid v\in V_1\setminus\{0\},\; v\perp
\Pi _{H_1}X\,\}\\ \label{eq:A.4}
&\ge  \min 
\{\,s_2(v,v)/\|v\|_{H_2}^2\mid v\in V_1\setminus\{0\},\; v\perp
                                   X\,\}\\ \nonumber
&\ge  \min \{\,s_2(v,v)/\|v\|_{H_2}^2\mid v\in V_2\setminus\{0\},\; v\perp X\,\},
\end{align} 
since $s_1(v,v)\ge s_2(v,v)$ on $V_1$, and $ V_2$ contains more
elements than $V_1$.
Taking the maximum over all  subspaces  $X$ of $H_2$ of dimension
$\le j-1$,
we get the $j$'th eigenvalues, which then must satisfy
the inequality \eqref{eq:A.2}.
\end{proof}

\section{Acknowledgement} The author thanks H.\ Schlichtkrull
for geometrical information, and thanks the referee for valuable criticism.



\end{document}